\documentclass[12pt]{amsart}

\usepackage{amsmath,amssymb,amsfonts,amsthm,amscd,indentfirst}

\usepackage{indentfirst}
\usepackage{hyperref}
\usepackage{graphicx}

\usepackage{color}
\usepackage[dvipsnames]{xcolor}

\usepackage{subcaption}

\numberwithin{equation}{section}

\usepackage{setspace}
\usepackage{amsmath, amssymb,amsthm, amsfonts, color}

\newtheorem{thm}{Theorem}[section]

\newtheorem{prop}{Proposition}[section]

\theoremstyle{definition}
\newtheorem{definition}{Definition}[section]

\theoremstyle{remark}
\newtheorem{remark}{Remark}[section]

\newcommand{\tr}{\mbox{tr}}

\renewcommand{\div}{\mbox{div}}

\newcommand{\R}{\mathbb R}

\newcommand{\be}{\begin{equation}}
\newcommand{\ee}{\end{equation}}
\newcommand{\bee}{\begin{equation*}}
\newcommand{\eee}{\end{equation*}}
\newcommand{\bal}{\begin{aligned}}
\newcommand{\eal}{\end{aligned}}

\def\p{\partial}

\def\la{\langle}
\def\ra{\rangle}

\def\Pi{\displaystyle{\mathbb{II}}}
\renewcommand\({\left(}

\def\vh{\vspace{.2cm}}

\def\m{\mathfrak{m}}

\def\F{\mathcal{F}}
\def\E{\mathcal{E}}

\def\hg{\bar g}

\def\CL{C^3_{L}}

\begin{document}

\title{Interpreting Mass via Riemannian Polyhedra}

\author{Pengzi Miao}
\address[Pengzi Miao]{Department of Mathematics, University of Miami, Coral Gables, FL 33146, USA}
\email{pengzim@math.miami.edu}

\thanks{P. Miao's research was partially supported by NSF grant DMS-1906423.}

\begin{abstract}
We give an account of some recent development that connects the concept of mass in general relativity 
to the geometry of large Riemannian polyhedra, in the setting of both asymptotically flat and asymptotically 
hyperbolic manifolds.
\end{abstract}

\maketitle

\markboth{Pengzi Miao}{Interpreting mass via Riemannian polyhedra}

\section{Scalar curvature and mass in general relativity}

The scalar curvature of a Riemannian manifold is a basic curvature quantity of the metric. 
Besides prolific motivations from Riemannian geometry,  
an important incentive to study the scalar curvature 
is from Einstein's theory of general relativity. 
We briefly outline this relativistic role of scalar curvature below. 

Suppose $ \mathcal{S}^4$ is a spacetime, i.e. a time-oriented Lorentz $4$-manifold. 
Suppose $ M^3$ is a spacelike hypersurface with a unit normal $n$ in $\mathcal{S}$. 
By the Einstein equation and the Gauss equation,
\be \label{eq-const}
R (g) - | K|_g^2 + ( \tr_g K )^2 = 16 \pi \, T (n, n) .
\ee
Here $ R(g)$ is the scalar curvature of  the Riemannian metric $g$ on $M$, 
$ K $ is the second fundamental form of $M$ in $\mathcal{S}$, and $T $ is the stress energy tensor 
determined by matter distribution in  $\mathcal{S}$.

A physically reasonable spacetime often satisfies suitable  energy conditions 
which guarantee  $T(n, n) \ge 0 $. In this case, 
if $M$  is a maximal slice, i.e. $ \tr_g K = 0 $, then \eqref{eq-const} implies  
$ R (g)   \ge 0 $. 
In general, under an assumption that  $\mathcal{S}$ satisfies the dominant energy condition, study of the triple $(M, g, K)$
often can be related to the case of $R(g) \ge 0$ (see \cite{SchoenYau81} for instance).

If the spacetime $\mathcal{S}$ models an isolated system, the Riemannian manifold $(M, g)$ is often assumed to be asymptotically flat. 
For such manifolds, a fundamental result on nonnegative scalar curvature metrics is the Riemannian positive mass theorem.

\begin{thm} [Schoen-Yau\cite{SchoenYau79}, Witten\cite{Witten81}] \label{thm-pmt}
Suppose $(M, g)$ is a complete, asymptotically flat $3$-manifold with nonnegative scalar
curvature. Then the mass $ \m (g) $ of $(M, g)$ satisfies $ \m (g) \ge 0 $, and equality holds 
 if and only   if $(M, g)$ is isometric to the Euclidean space. 
\end{thm}
 
We recall the definitions of asymptotically flat manifolds and their mass.

\begin{definition} \label{def-AF}
A Riemannian manifold $(M^n, g) $ is said to be asymptotically flat (with one end)
if there exists a compact set $K $ such that $ M \setminus K $ is diffeomorphic to 
$ \R^n \setminus B_r (0) $, where $ B_r (0)  = \{ |x| < r \} $ for some $ r > 0 $, and
with respect to the standard coordinates $\{ x_i \}$ on $\R^n$, 
the metric $g$ satisfies $ g_{ij} = \delta_{ij} + h_{ij}$, where 
\be \label{eq-g-decay}
h_{ij} =  O ( | x |^{-p} ) , \ \p h_{ij}  =  O ( | x |^{ - p - 1}) , \  
\p \p h_{ij}  =  O ( | x |^{ - p - 2 } )
\ee 
for some $ p > \frac{n-2}{2} $.
Moreover, the scalar curvature of $g$ is integrable on $(M, g)$.
\end{definition}

For an asymptotically flat $ (M^n, g) $, its mass 
(introduced by Arnowitt, Deser and Misner \cite{ADM}), 
is defined as
\be \label{eq-mass-def}
\m (g) = \lim_{ r \rightarrow \infty} \frac{1}{ 2 (n-1) \omega_{n-1} } \int_{ S_r }  (g_{ij,i} - g_{ii,j} ) \, \bar \nu^j \, d \bar \sigma .
\ee
Here $ S_r $ is the coordinate sphere $ \{  | x | = r \}$, 
$ d \bar \sigma$ denotes the Euclidean volume element on $ S_r$, 
$ \bar \nu $ is the Euclidean  outward unit normal to $ S_r$, 
and summation is applied over repeated indices. 
It was proved by Bartnik \cite{Bartnik86}, and also by Chru\'sciel \cite{Chrusciel},  that $\m(g)$
is independent on the choice of coordinates satisfying \eqref{eq-g-decay}.

For dimensions $ n \le 7$, Schoen had a proof of the positive mass theorem in \cite{Schoen-LNM}.
For arbitrary dimensions, the theorem was shown by Schoen-Yau \cite{SchoenYau17}.
A different approach was provided by Lohkamp \cite{Lo1, Lo2, Lo3, Lo4}.

\section{Gromov's geometric comparison theory for scalar curvature}

In Riemannian geometry, comparison theorems are useful tools in understanding manifolds with 
various notions of lower curvature bounds. 
In two dimensions, scalar curvature reduces to 
Gauss curvature, and geometric triangles were compared to those in a model space
(see \cite{Aleksandrov48, Aleksandrov51, Toponogov} for instance).

In \cite{Gromov14}, Gromov suggested Riemannian polyhedra play a useful role 
 in the study of scalar curvature  in higher dimensions. 
Given a convex polyhedron $P$ in  $ (\R^n, \bar g)$, the Euclidean space, 
Gromov \cite{Gromov14, Gromov18a} conjectured that,  if $g$ is a Riemannian 
metric on $P$ satisfying
\begin{itemize}
\item  $g$ has nonnegative scalar curvature;
\item each face of $ P$ has nonnegative mean curvature in $(P, g)$; and
\item along each edge of $ P$, the dihedral angle of $(P, g)$ is 
less than or equal to the corresponding dihedral angle of $P$ in $(\R^n, \bar g)$,
\end{itemize}
then $(P, g)$ is isometric to a polyhedron in $ (\R^n, \bar g)$.

\begin{figure}[h]
\centering
\includegraphics[scale=.7]{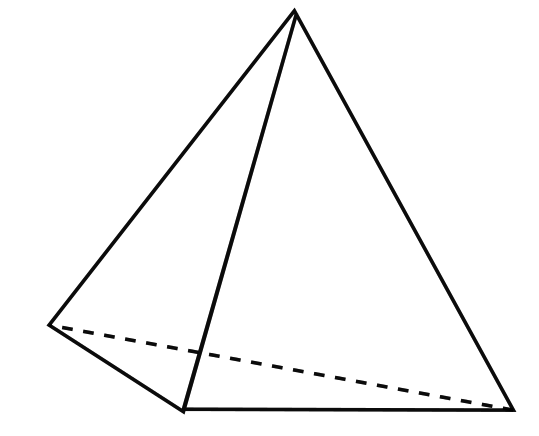}
\caption{A Riemannian tetrahedron}
\label{fig-t}
\end{figure}

Significant progress toward this conjecture was obtained by Li \cite{Li17, Li20, Li20h}.
In $3$-dimension, Li \cite{Li17} proved the conjecture for a wide class of polyhedra including cubes and 
$3$-simplices. In higher dimensions, Li \cite{Li20} showed the conjecture holds for a class of prisms with dimensions less than or equal to $7$. 
In \cite{Li20h}, Li extended the comparison theory to metrics with  a negative scalar curvature lower bound.

An immediate connection between Gromov's conjecture and the positive mass theorem is that
the conjecture being valid for a single polyhedron implies the latter.
This is because, by the work of Lohkmap \cite{Lo99}, if the mass of an asymptotically flat $(M, g)$ is negative, 
the metric $g$ can be deformed to another metric $\tilde g$ with nonnegative scalar curvature, which is 
exactly Euclidean outside a compact set. 
In this case, by considering a large polyhedron enclosing the nontrivial geometry 
in the compact set, one would obtain a contradiction to Gromov's conjecture
(see \cite[Theorem 5.2]{Li20} for details).

\section{Mass and coordinate cubes in $3$-dimension}

For an asymptotically flat $3$-manifold $(M, g)$, a connection between the mass $\m(g)$ and 
large Riemannian cubes was found in \cite{Miao19}.

\begin{thm}[\cite{Miao19}] \label{thm-cube-mass}
Let $ (M^3, g)$ be an asymptotically flat $3$-manifold. Given any large  constant $ L > 0 $, 
let $ \CL$ denote a coordinate cube with side length $2L $ centered at the coordinate origin.
Let $\mathcal{F} $ and $ \mathcal{E}$ denote the sets of the faces and the edges of $ \CL$, 
respectively.  Then, as $L \to \infty$, 
\be \label{eq-Gromov-mass-1}
\m (g) =  -\frac{1}{8\pi} \int_{\mathcal{F} } H \, d \sigma + \frac{1}{8\pi} \int_{\mathcal{E}} \left( \alpha  -  \frac{\pi}{2}  \right) \, d s + o(1).
\ee
Here $H$ is the mean curvature of each face with respect to the outward unit normal, 
$ \alpha $ is  the dihedral angle between the two adjacent faces  along each edge, 
$ d \sigma$ and $ d s$ are the area and the length measure with respect to $g$, respectively. 
\end{thm}

If  $(M^3, g)$ is complete with nonnegative scalar curvature,
the positive mass theorem (Theorem \ref{thm-pmt})  shows $ \m (g) > 0 $, unless 
$ (M, g)$ is isometric to $ \R^3$. This combined with \eqref{eq-Gromov-mass-1} shows
\be \label{eq-improve-Gromov}
-\frac{1}{8\pi} \int_{\mathcal{F}} H \, d \sigma + \frac{1}{8\pi} \int_{\mathcal{E}} \left( \alpha - \frac{\pi}{2} \right) \, d s \ge 0, 
\ \text{for large} \ L.
\ee

A related but yet different formula,  relating $\m(g)$ to $C^3_L$,
was also derived in \cite{Miao19}.

\begin{thm}  [\cite{Miao19}] \label{thm-GB}
Let $(M^3, g)$, $ L$ and $ \CL$ be given as above. 
For each $ t \in [- L, L]$ and each $ k \in \{1, 2, 3\}$, 
let $S_t^{(k)}$ be the curve given by the intersection between $\partial C^3_L$ 
and the coordinate plane  $\lbrace x_k=t\rbrace$.
Then, as $ L \to \infty$,  the mass of $(M^3,g)$ satisfies 
\begin{equation}\label{eq-GB}
\m (g) = \frac{1}{8\pi}\sum_{k=1}^3 \int_{-L}^L \m^{(2)}_k  ( t, L) \, dt + o (1) .
\end{equation}
Here
\be
\m^{(2)}_k (t, L) =  2\pi -\int_{S_t^{(k)}}\kappa^{(k)}ds -\beta^{(k)}_t  , 
\ee
$\kappa^{(k)}$ is the geodesic curvature of $S_t^{(k)}$ in $\lbrace x_k=t\rbrace$, 
and $\beta^{(k)}_t$ is the sum of the turning angle of $S_t^{(k)}$ at its four vertices. 
\end{thm}

The quantity $ \m^{(2)}_k (t, L) $ measures the angle deficit of the (large) portion of 
the coordinate plane $\{ x^k = t \}$ inside the cube $\CL$.
\eqref{eq-GB} shows the mass of $(M^3, g)$  can be obtained by suitably integrating 
such angle deficit associated to all coordinate planes.
It is worthy of noting that, in the setting of asymptotically conical surfaces, the angle deficit can be used 
to define the $2$-dimensional ``mass" of those surfaces (see \cite{Chr, CL2019} for instance).   

Formula \eqref{eq-GB} was largely prompted by a recent monumental work of Stern \cite{Stern19} 
that relates the scalar curvature and the level sets of harmonic functions (and maps). 
Together with the Gauss-Bonnet theorem and Stern's formula from \cite{Stern19}, \eqref{eq-GB}
provides insight into the recent new proof of the $3$-dimensional positive mass theorem by 
Bray-Kazaras-Khuri-Stern \cite{BKKS19}. 
Below we briefly outline this implication following \cite[Remark 5]{Miao19}.

Stern's formula \cite[equation (14)]{Stern19}, in its simplest form, shows
\be \label{eq-Stern-1}
 \Delta | \nabla u | =   \frac{1}{2   | \nabla u | }   \left[     | \nabla^2 u |^2  +    | \nabla u |^2 ( R - 2 K_{_\Sigma} )  \right],
\ee
near points where $ \nabla u \ne 0$. Here 
 $u$ is a harmonic function on a $3$-manifold $(M^3, g)$, $ R$ and $K_{_\Sigma}$ denote
the scalar curvature of $g$  and the Gauss curvature of $\Sigma$, the level set of $u$, respectively. 

Suppose  the manifold $(M^3, g)$ is asymptotically flat. An observation of Bartnik \cite{Bartnik86} was
\be \label{eq-Bartnik-formula}
 \sum_{i=1}^3 \int_{S_\infty}    \frac12 \frac{\p }{\p \nu} | \nabla y_i |^2 \, d \sigma 
  = 16 \pi \m (g) .
\ee
Here $\{ y_i \}$ are harmonic coordinates on $(M, g)$ near infinity,
and $ \int_{S_\infty}$ represents the limit of integration along a sequence of suitable surfaces  approaching infinity.
Since $ | \nabla y_i | $ tends to $ 1$ suitably fast, it is easily checked that
\eqref{eq-Bartnik-formula} is equivalent to 
\be \label{eq-Bartnik-formula-1}
 \sum_{i=1}^3  \int_{S_\infty}  \frac{\p }{\p \nu} | \nabla y_i | \, d \sigma = 16 \pi \m (g).
\ee

Considering \eqref{eq-Stern-1} and \eqref{eq-Bartnik-formula-1}, 
one may want to seek a formula that computes $\m (g)$ in terms of the  level sets of $y_i$.
This was the motivation to the derivation of \eqref{eq-GB} in \cite{Miao19}. 
It turned out that \eqref{eq-GB} itself does not require $\{ x_i \}$ to be harmonic. 

Next,  for simplicity, suppose $M$ has no boundary. Consider a harmonic map 
\be
U = (u_1, u_2, u_3): (M, g) \rightarrow (\R^3, \bar g ) 
\ee
such that $U$ is a diffeomorphism near infinity, 
and, in this coordinate chart defined by $U$, $ g$ satisfies \eqref{eq-g-decay}. 
By the usual construction of harmonic coordinates (for instance Theorem 3.1 in \cite{Bartnik86}), such a map $U$ always exists.
If one further supposes the regular level set $\Sigma_t^{(i)}$ of each $u_i$ satisfies 
$\chi (\Sigma_t^{(i)}) \le 1 $ (which is the case if $ M $ is $ \R^3$ for instance),  
then  it follows from Stern's formula \eqref{eq-Stern-1},  Bartnik's formula \eqref{eq-Bartnik-formula-1}, 
the mass formula \eqref{eq-GB}, and the Gauss-Bonnet theorem that
\be \label{eq-PMT-refined}
\begin{split}
& \  24 \pi \, \m (g) \\
= & \  \lim_{L \to \infty} \sum_{k=1}^3 \int_{\p \CL} \frac{\p}{\p \nu} | \nabla u_k | \, d \sigma  +  
 \lim_{L \to \infty}  \sum_{k=1}^3  \int_{-L}^L  \left[    ( 2 \pi - \beta^{ \la k \ra }_t ) - \int_{ C^{ (k)}_t   } \kappa^{(k)}  \, d s  \right]  \, dt  \\
 \ge & \ \sum_{i=1}^3 \ \int_{M} \frac{1}{2}  \left[ \frac{1}{  | \nabla u_i | }  | \nabla^2 u_i |^2 + R |\nabla u_i | \right] d V.
 \end{split}
\ee
In particular, if $ R \ge 0 $, then $ \m(g) \ge 0 $.
We remark that inequality \eqref{eq-PMT-refined} is weaker  than the result in \cite{BKKS19}. 
The lower bound of $\m(g)$ in Theorem 1.2 of \cite{BKKS19} 
 only used a single harmonic function. 

A potentially useful feature of the above mapping setting  \eqref{eq-PMT-refined} is the following.
Suppose $ g_{ij} = g (\p_{u_i}, \p_{u_j} ) $. 
Given any constant $ \Lambda > 0 $, let  $ r_o > 0 $ be a constant such that 
$$ | g_{ij}  |  \le \Lambda \ \ \mathrm{and} \ \ | g^{ij} | \le \Lambda ,   \ \ \mathrm{if} \
 | u | \ge r_o  . $$
Using the fact 
$ | \nabla u_i |^2 = g^{ii} $ and  $ \nabla^2 u_i (\p_{u_k}, \p_{u_l} ) = - \Gamma^i_{kl} $,
we see from \eqref{eq-PMT-refined} that, if the scalar curvature is nonnegative, then 
\be
\m (g) \ge C \int_{ | u | \ge r_o} \sum_{i,j,k} | \p_i g_{jk} |^2 ,
\ee 
where $C$ is a constant depending only on $\Lambda$. 
This shows, outside a large $u$-coordinate sphere $\{ |u| = r_o \}$, 
the $L^2$-norm of $\p g$ is controlled by the mass. 

\section{Mass and polyhedra in general dimensions} \label{sec-mass-poly}

Motivated by Theorem \ref{thm-cube-mass}, Piubello and the author \cite{MiaoPiubello21} 
considered a sequence of polyhedra of general types, approaching infinity in an asymptotically flat manifold 
of arbitrary dimensions. 

\begin{figure}[h]
\centering
\includegraphics[scale=.35]{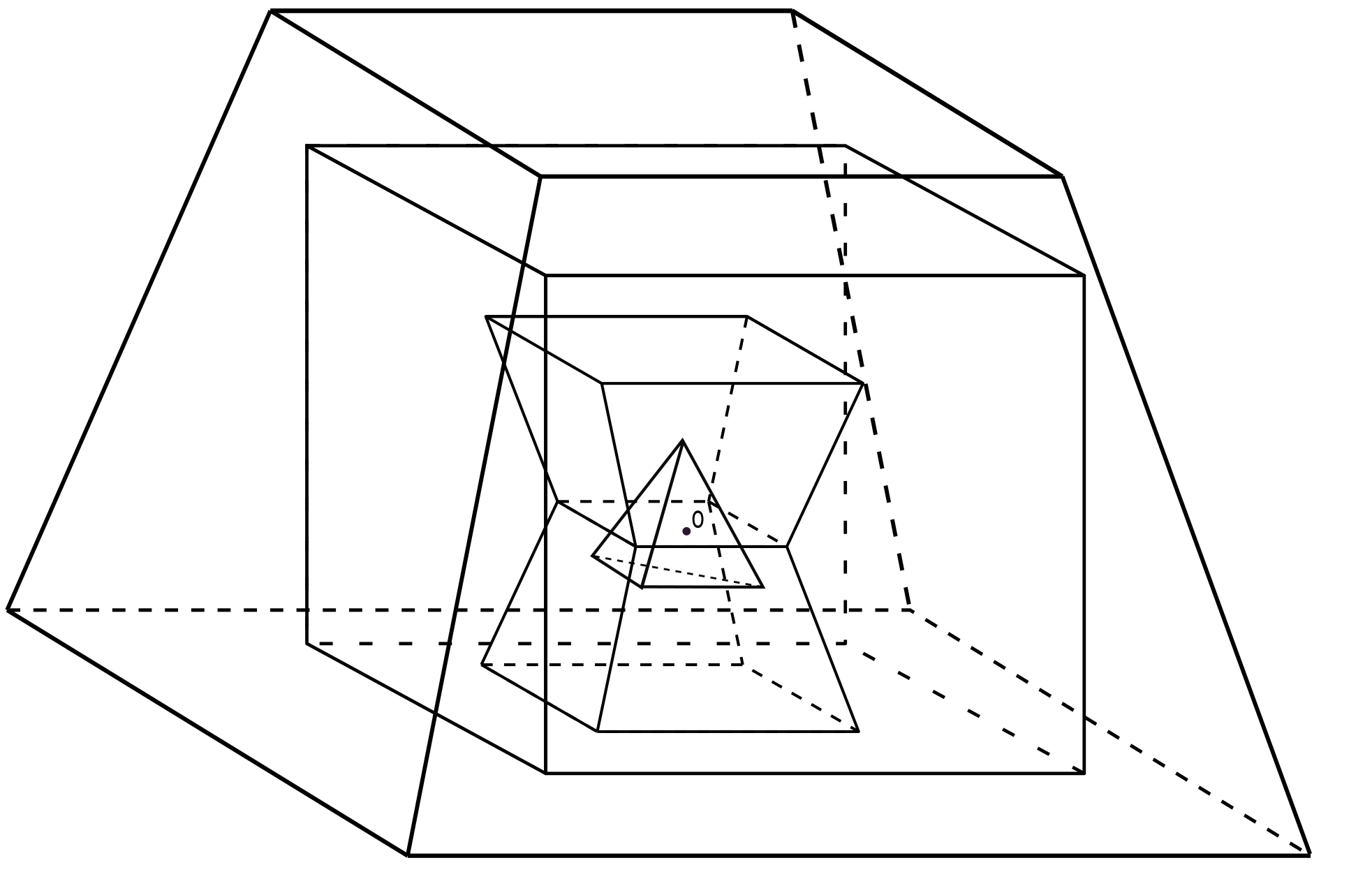}
\caption{A sequence of polyhedra approaching $\infty$ in $(M^n, g)$}
\label{fig-p}
\end{figure}

\begin{thm} [\cite{MiaoPiubello21}]  \label{thm-poly-mass}
Let $(M^n, g)$ be an asymptotically flat manifold with dimension $ n \ge 3$. 
Let $ \{ P_k \}$ denote a sequence of Euclidean polyhedra in a coordinate chart $\{ x_i \}$ 
that defines the end of $(M, g)$. Suppose $ \{ P_k \}$ satisfies the following conditions:
\begin{itemize}
\item[a)]  $ P_k $ encloses the coordinate origin  and $ \lim_{k \to \infty} r_{_{P_k}} = \infty $, where
$$ r_{_{P_k} } = \min_{x \in \p P_k}  | x |  ; $$
\item[b)] $ | \F (\p P_k)  | = O ( {r^{n-1}_{_{P_k}}} )$, where $\F (\p P_k) $ denotes  the union of all the faces in $\p P_k$, 
and $ | \F ( \p P_k)  | $ is the Euclidean $ (n-1)$-dimensional volume of $ \F (\p P_k) $;
\vh
\item[c)] $ | \E  ( \p P_k)  |  = O ( r^{n-2}_{_{P_k} })$, where $ \E ( \p P_k) $ denotes the union of all the edges in $ \p P_k$, 
and  $| \E ( \p P_k)  | $ is  the Euclidean $(n-2)$-dimensional volume of $ \E ( \p P_k) $; and
\vh
\item[d)] along each edge, the Euclidean dihedral angles $\bar \alpha $ of $P_k$ satisfies 
$$  | \sin \bar \alpha |   \ge c $$ 
for some constant $ c > 0$ that is independent on $k$.
\end{itemize}
Then the mass of $(M, g)$ satisfies, as $ k \to \infty$,
\be \label{eq-main}
 \m (g)  = \frac{1}{(n-1) \omega_{n-1} } \left( - \int_{\F (\p P_k) } H \, d \sigma + \int_{\E (\p P_k) } ( \alpha - \bar \alpha) \, d \mu \right)  + o(1) . 
\ee
Here $ \omega_{n-1}$ is the volume of the standard $(n-1)$-dimensional  sphere, 
 $ H$ denotes the mean curvature of the faces $\F (\p P_k) $ in $(P_k, g)$, $\alpha$ denotes the dihedral angle along 
the edges $\E (\p P_k) $ of $(P_k, g)$, and $ d \sigma$, $ d \mu $ denote the volume element on $\F (\p P_k) $, $ \E ( \p P_k) $, respectively, with respect to the induced metric from $g$.
\end{thm}

The sequence $ \{ P_k \} $ in Theorem \ref{thm-poly-mass} has the feature that 
it can consist of different types of polyhedra.
Moreover, elements in $\{ P_k \}$ can be non-convex. 
The angle condition 
$$ | \sin \bar \alpha | \ge c  $$
is imposed due to $ y' (x)  = - \frac{1}{\sin y} $ if $ y = \arccos x $. 
In the proof of Theorem \ref{thm-poly-mass} in \cite{MiaoPiubello21}, 
one obtains the dihedral angle from the metric coefficients. A condition of $ | \sin \bar \alpha | \ge c $ 
helps to convert estimates on the metric to estimates on the dihedral angle. 

\begin{remark}
Let $P$ be a fixed polyhedron in $ \R^n$ that encloses the coordinate origin.
Let $ P_{(r)} $ be the polyhedron obtained by scaling $P$ with a constant factor $r $ in the coordinate space. 
For large $r$, the family $\{ P_{(r)} \}$ 
satisfies all the conditions a) -- d) in Theorem \ref{thm-poly-mass}.  As a result,  \eqref{eq-main} 
is applicable to such a family  $\{P_{(r)} \}$ as $ r \to \infty$. 
\end{remark}

\begin{remark} \label{rem-angle-AF}
Formula \eqref{eq-main} implies the limit of 
$ - \int_{\F (\p P_k) } H \, d \sigma + \int_{\E (\p P_k) } ( \alpha - \bar \alpha) \, d \mu  $
is geometric, i.e. independent on the choices of $\{ x_i \}$ and $ \{ P_k \}$.  It is not hard to see 
this is not the case for the limit of each individual factor. For instance, let $(M^n, g_m)$ denote the spatial 
Schwarzschild manifold. If one uses the isotropic coordinates $\{ x_i \}$, in which 
$ g_m = \left( 1 + \frac{m}{2 |x|^{n-2} }  \right)^\frac{4}{n-2} \bar g , $
and let $ \{ P_k \}$ be the associated coordinate cubes, then 
$ \alpha = \bar \alpha  $ due to the fact that $g$ and $\bar g$ are conformal  
in these $\{x_i\}$ coordinates. In particular,  $\int_{\E (\p P_k) } ( \alpha - \bar \alpha) \, d \mu = 0 $
in this case.
On the other hand, if one uses the spherical coordinates $(r, \theta)$, in which 
$ g_m = \frac{1}{ 1 - \frac{2m}{ r^{n-2} } } d r^2 + r^2 \sigma_o $, where $ \sigma_o$ is the standard  
metric on the unit sphere $\mathbb{S}^{n-1}$ and $ \theta \in \mathbb{S}^{n-1}$, then 
$ g$ is not conformal to the background Euclidean metric $\bar g$ in the rectangular coordinates $\{ y_i \}$
given  by $ y  = r \theta $, unless $m =0$. 
Suppose $ m \ne 0 $ and choose  $ \{ P_k \}$ as the coordinate 
cubes in these $\{ y_i \}$ coordinates. One can easily check that, passing to the limit,
the angle deficit term $\int_{\E (\p P_k) } ( \alpha - \bar \alpha) \, d \mu $ has a nontrivial contribution to the mass.
\end{remark}

We briefly comment on the proof of Theorem \ref{thm-poly-mass}. 
Formula \eqref{eq-main} was essentially revealed by a special structure that is present 
in the linearization of the mean curvature function. 

To explain this structure, we fix some notations. 
Let $ \Sigma^{n-1} \subset M^n$ be a given hypersurface.
Let $g$, $ \bar g$ be two Riemannian metrics on $M$, and 
let $\gamma$, $\bar \gamma$ denote the metrics on $ \Sigma$ induced from $g$, $\bar g$, respectively. 
Let $\nu$ and $ \bar \nu$ be the corresponding unit normal vectors to $ \Sigma$ in $(M, g)$ and $(M, \bar g)$. 
We assume $\nu$ and $\bar \nu $ point to the same side of $\Sigma$. 
Let $ A$, $ H$ and $ \bar A$, $ \bar H $ denote the second fundamental form, 
the mean curvature function of $\Sigma$ in $(M, g)$ and $(M, \bar g)$, respectively.  
Viewing $\bar g$ as a background metric, we write $ g = \bar g + h $ and assume
$   | h |_{\bar g} < \epsilon (n) $, 
where $ \epsilon (n)  $ is a small positive constant depending only on $n$.
Here $ | \cdot |_{\bar g}$ denotes the norm with respect to $\bar g$.

\begin{prop}[\cite{MiaoPiubello21}]  \label{prop-mean-curvature}
The mean curvature functions $H$ and $ \bar H$ satisfy 
\be \label{eq-H-bar-H-prop}
\begin{split}
2 ( H - \bar H ) 
= & \ ( d \, \mathrm{tr}_{\bar g} \, h - \mathrm{div}_{\bar g} \, h ) (\bar \nu) - \mathrm{div}_{\bar \gamma} X - \la h, \bar A \ra_{\bar \gamma}  \\
& \  + | \bar A |_{\bar g} \, O ( | h |^2_{\bar g} ) + O ( | \bar D h |_{\bar g} \, | h |_{\bar g} ).
\end{split}
\ee
Here $ \mathrm{div}_{\bar g} (\cdot)$ and $ \mathrm{tr}_{\bar g} (\cdot) $ denote the divergence and trace on $(U, \bar g)$;
$ \mathrm{div}_{\bar \gamma} (\cdot) $ denotes the divergence on $(\Sigma, \bar \gamma)$;
 $X $ is the vector field on $ \Sigma$ that is dual to the $1$-form $h(\bar \nu, \cdot) $ with respect to  $\bar \gamma$;
 and $\bar D$ denotes the covariant derivative with respect to $\bar g$.
Given two functions $f$ and $\phi$, we write $ f =  O ( \phi )$ to denote $ | f  | \le C | \phi | $ with a constant $ C $ that depends only on $n$.
\end{prop}

We see that, on the right side of \eqref{eq-H-bar-H-prop},  
the first term $( d \, \tr_{\bar g} \, h - \div_{\bar g} \, h ) (\bar \nu) $ 
relates to the mass integral; the divergence term $ \div_{\bar \gamma} X$, upon integration by parts, produces 
the angle along the edges; and the third term, involving the second fundamental form $\bar A$, 
indicates the advantage of making use of a polyhedron, whose faces are totally geodesic. 

Similar to \eqref{eq-improve-Gromov}, 
if  $(M^n, g)$ is complete, asymptotically flat and has nonnegative scalar curvature, 
it follows from the positive mass theorem and \eqref{eq-main} that
\be \label{eq-main-cor}
- \int_{\F (\p P_k) } H \, d \sigma + \int_{\E (\p P_k) } ( \alpha - \bar \alpha) \, d \mu \ge 0 , \ \text{for large } k .
\ee
Here $ \{ P_k \}$ is any sequence of Euclidean polyhedra satisfying conditions in Theorem \ref{thm-poly-mass}.
In particular, for a fixed Euclidean polyhedron  $ P$, 
\be
- \int_{\F_{(r) }} H \, d \sigma + \int_{\E_{(r) }} ( \alpha - \bar \alpha) \, d \mu \ge 0 , \ \text{for large } r , 
\ee
where $\F_{(r)}$ and $ \E_{(r)}$ are the faces and edges of the polyhedron $P_{(r)}$ that is obtained by 
scaling $P$ by a constant $r$.

As illustrated by Theorem \ref{thm-GB}, it is worthy of  taking $ \{ P_k \}$ 
as a sequence of large coordinate cubes in $(M^n, g)$.
Cubes have a feature that, when sliced by hyperplanes parallel to their face, 
the resulting intersections are large $(n-1)$-dimensional cubes.

\begin{figure}[ht]
\centering
\includegraphics[scale=.5]{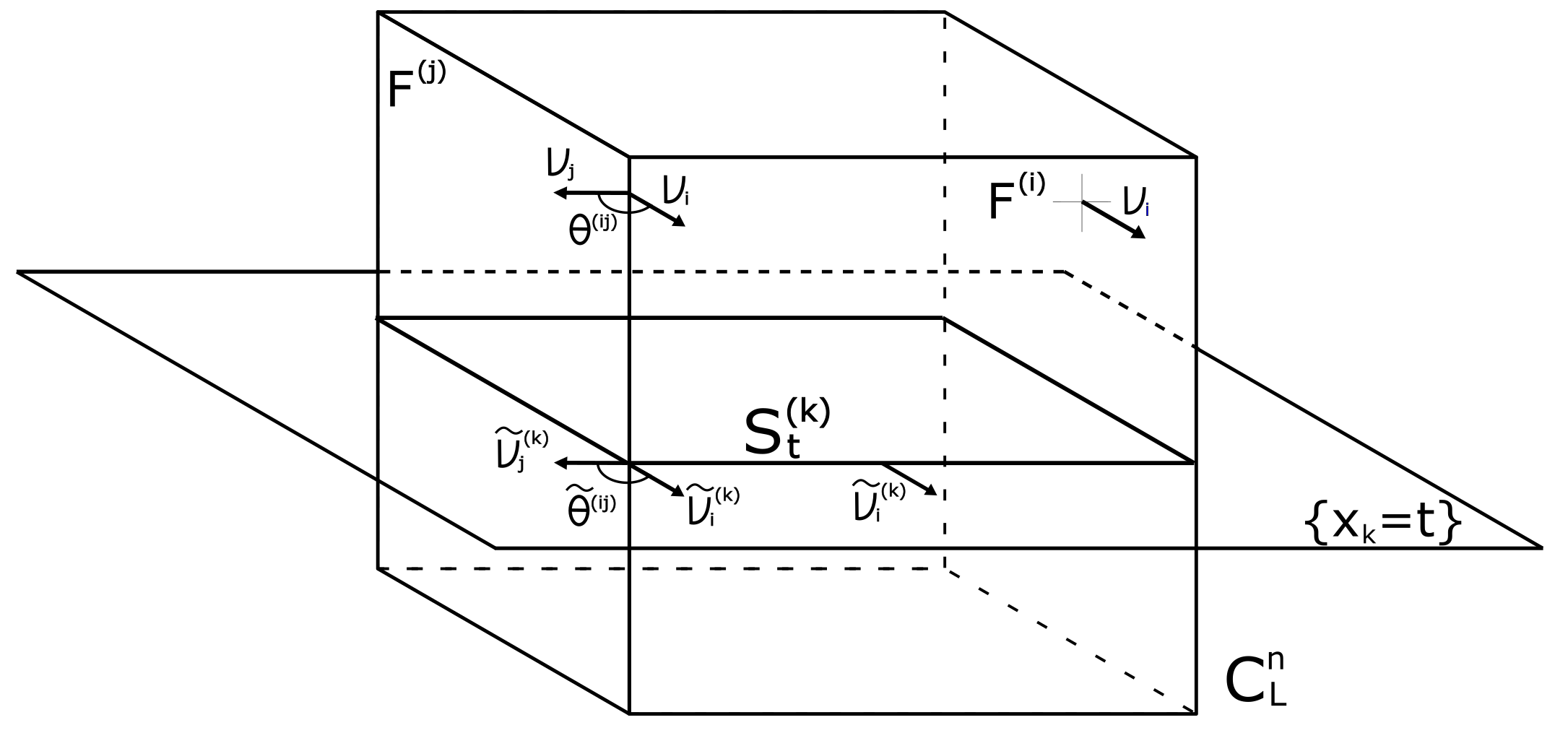}
\caption{Slicing a coordinate cube}
\label{fig-cube}
\end{figure}

The following can be viewed as a higher dimensional analog of Theorem \ref{thm-GB}. 

\begin{thm} [\cite{MiaoPiubello21}] \label{thm-slicing}
Let $(M^n, g)$ be an asymptotically flat manifold with dimension $ n \ge 4$.
Given any $k \in \{ 1, \cdots, n \}$,  any large constant $L$, and any $ t \in [-L, L]$,  
there is a quantity $ \m_k^{(n-1)} (t, L)$, associated to the coordinate hyperplane 
$\{ x_k = t \}$, defined in \eqref{eq-def-mktl} below, such that 
\be
\m ( g|_{  \{ x_k = t \} } ) = \lim_{ L \to \infty } \m_k^{(n-1)} (t, L) ,
\ee
and
\begin{equation}\label{eq-slicing}
  \m (g)=  \frac{ \omega_{n-2} }{ (n-1) \omega_{n-1} }   \lim_{ L \to \infty} \sum_{k=1}^n \int_{-L}^L \m^{(n-1)}_k ( t, L ) \, dt .
\end{equation}
\end{thm}

We explain the quantity $\m_k^{(n-1)}( t, L) $.  
Let $C^n_L$ denote a coordinate cube in $(M, g)$, 
with its center at the coordinate origin and side (coordinate) length $2L$. 
For every $ t  \in [-L, L]$ and $k \in \{ 1, \cdots, n \}$, 
let $S_t^{(k)}$ be the intersection of $\partial C^n_L$ with the coordinate hyperplane $\{ x_k = t \}$, i.e.
$$S_t^{(k)}=\partial C^n_L\cap \lbrace x_k=t\rbrace.$$
Then $S_t^{(k)}$ is the boundary of an $(n-1)$-dimensional  cube in $\{ x_k = t \}$ (see Figure \ref{fig-cube}).
On the hypersurface $\{ x_k = t \}$,  let  $\tilde{H}^{(k)}$ be the mean curvature of 
each face of $S_t^{(k)}$ with respect to the outward normal $ \tilde \nu^{(k)}_i$ and 
let $\tilde{\alpha}^{(k)}$ be the dihedral angle along every edge of $ S^{(k)}_t$, 
both taken with respect to $g$.

Associated to each $ S^{(k)}_t$ in $\{ x_k = t \}$, $ m_k^{(n-1)} (t, L)  $ is given by 
\be \label{eq-def-mktl}
\begin{split}
& \ \m_k^{(n-1)} (t, L) \\
=  & \  \frac{1}{ (n- 2)\omega_{n-2}  } \left(  - \int_{ \F ( S^{(k)}_t ) } \tilde{H}^{(k)} d\sigma^{n-2}  + \int_{\mathcal{E}(S_t^{(k)})}\( \tilde\alpha^{(k)} - \frac{\pi}{2}\right)d\mu^{n-3} \right) ,
\end{split}
\ee
where  $ d \sigma^{n-2}$ and $ d \mu^{n-3}$  denote the volume forms, 
induced from the metric $g$, on the faces and edge of $S^{(k)}_t$.
The assertion 
$$ \m ( g|_{  \{ x_k = t \} } ) = \lim_{ L \to \infty } \m_k^{(n-1)} (t, L) $$
is a direct consequence of  Theorem \ref{thm-poly-mass},
where $ g |_{\{ x_k = t \} }$ denotes the metric on $\{ x_k = t \}$ induced from $g$.
In many cases, for instance if $g$ has a decay rate of $p > n - 3$, this limit will be zero 
as $ \m ( g|_{  \{ x_k = t \} } ) = 0 $.

\vspace{.2cm}

A couple of questions arise in relation to Theorems \ref{thm-poly-mass} and \ref{thm-slicing}.

\begin{itemize}

\item[a)]  The angle condition $ | \sin \bar \alpha | \ge c $ serves only as a sufficient condition for \eqref{eq-main}.
It would be good to weaken this condition. For instance, a replacement of $ | \sin \bar \alpha | \ge c $ by a suitable slow 
decay condition of $ \sin \bar \alpha $, relative to the growth of the volume of the faces and edges, may be desirable. 

\vspace{.2cm}

\item[b)] In dimension three,  formula \eqref{eq-GB} combined with the Gauss-Bonnet theorem and 
Stern's work [17] yielded a proof of the positive mass theorem.  
Could formula \eqref{eq-slicing} play a role in any new proof of the positive mass theorem 
in high dimensions?

\end{itemize}

\section{Mass of asymptotically hyperbolic manifolds} 

The preceding connection between polyhedra and the mass of asymptotically flat manifolds prompts one to consider
if an analogous connection exists on manifolds that are asymptotically hyperbolic. 

We recall the definition of asymptotically hyperbolic manifolds and their mass functional 
following the work of Chru\'{s}ciel and Herzlich 
\cite{CH03}.
Let $ (\mathbb{H}^n , \bar g) $ denote the standard hyperbolic space with constant sectional curvature $-1$.
One model of $ \mathbb{H}^n$ is the upper hyperboloid in the Minkowski spacetime
$$ \mathbb{H}^n=\{(z,t)\in\mathbb{R}^{n,1} \ | \ z_1^2+\cdots+z_n^2-t^2=-1, t>0\}. $$
Let $r= | z | $.  
A Riemannian manifold $(M^n,g)$ is called asymptotically hyperbolic if there exist a compact set $K\subset M$ and 
a diffeomorphism $\Phi:M\setminus K\to \mathbb{H}^n\setminus D $, where $D $ is some compact set in $ \mathbb{H}^n$, 
such that $h : =(\Phi^{-1})^* g - \hg $ satisfies 
$$ |h|_{\hg} +| \bar D h|_{\hg}+| {\bar D}^2 h|_{\hg}=O(r^{-q}), \ \text{for some rate} \ q>\frac{n}{2}, $$
as $ r \to \infty$, and 
$$\int_M r \left[ R_g+n(n-1) \right] \, d V_g <\infty . $$ 
Here $ R_g $ and $d V_g$  are  the scalar curvature and the volume element of $g$, respectively. 
The mass functional of $(M, g)$ is a linear functional
\be
\m ( \cdot ) : \mathcal{S} (\mathbb{H}^n )   \longrightarrow \R ,
\ee
where  $\mathcal{S} (\mathbb{H}^n )$ is  the space of static potentials on $\mathbb{H}^n$, given by 
\be
\mathcal{S} (\mathbb{H}^n ) = \text{span} \{ t, z_1, \cdots, z_n \}, \ \text{where} \ t = \sqrt{ 1 + r^2},  
\ee
and,  $\forall \, V \in  \mathcal{S} (\mathbb{H}^n )  $, 
\begin{equation}\label{eqn:AH mass integral}
	\m (V)=\lim_{r\to\infty}\int_{ |z| = r } 
	\left[ V \, ( \div_{\hg} h- d \mathrm{tr}_{\hg} h ) + (\mathrm{tr}_{\hg} h)d V 
	-h({\bar \nabla} V,\cdot) \right] (\bar \nu)\,d \bar \sigma  .
\end{equation}
Here $  d \bar \sigma $ is the volume element on the coordinate sphere $\{ |z| = r \}$ in $\mathbb{H}^n $ 
and $\bar \nabla $ denotes the gradient on $(\mathbb{H}^n, \bar g)$.
Similar to the asymptotically flat case,  the limit of the above integral can be computed along other suitable exhaustion sequences of $M$ (see \cite{Michel11} for instance).

It is common to represent $\m (\cdot) $ by an $(n+1)$-vector 
$(p_0,p_1,\ldots,p_n)$ with 
$$
 p_0= \m ( t ), \ \ p_i= \m (z_i), \ \forall \, i=1,\ldots,n. 
$$
The Riemannian positive mass theorem for asymptotically hyperbolic manifolds states that, if 
 $R_g\ge -n(n-1)$, then 
$  p_0^2 \ge \sum_{i=1}^n p_i^2 $, 
and equality holds if and only if $(M,g)$ is isometric to $\mathbb{H}^n$. 
The proof of this theorem has a rich history. 
We refer readers to \cite{XW01, Zhang04, CH03, ACG08, Chrusciel:2018co, Chrusciel.2019, Huang20, Sakovich20} 
and references therein.

Given a function $V$,  let $\mathbb{U} (V) (\cdot) $ be the associated $1$-form 
\be
\mathbb{U} (V) (\cdot) =  
V \, ( \div_{\hg} h- d \tr_{\hg} h ) + (\tr_{\hg} h)d V -h({\bar \nabla} V,\cdot) .
\ee
The following formula, which can be viewed as a weighted version of Proposition \ref{prop-mean-curvature},
was derived by Jang and the author in \cite{JangMiao20}.

\begin{prop} [\cite{JangMiao20}] \label{prop-mass-1-form}
Given two metrics $ \hg $ and $g$, a function $V$, and a hypersurface $\Sigma \subset M $,  
let $\bar \nu,\nu$ denote the unit normal vectors to $\Sigma$ pointing the same side in $(M, \hg), (M,g)$, respectively. Then, for $|h|_{\hg}$ small, 
\begin{equation} \label{eq-weighted-H}
\begin{split}	
\mathbb{U}(V)(\bar \nu) = & 2 V (\bar H -H) -\mathrm{div}_{\bar \gamma} (VX)  
+ (\mathrm{tr}_{\bar \gamma} h)dV( \bar \nu )-V \langle \bar A,h\rangle_{\bar \gamma} \\
& \ +  V \left[  |\bar A|_{\hg}O(|h|_{\hg}^2)+O(|\bar D h|_{\hg} |h|_{\hg}) \right] .			
\end{split}
\end{equation}
Here $\bar H$, $H$, $\bar \gamma$, $\bar A$ and $X$ are defined as in Proposition \ref{prop-mean-curvature}.
\end{prop}

Formula \eqref{eq-weighted-H} suggests that when evaluating the total flux of $ \mathbb{U}(V) (\cdot) $ across
a hypersurface $ \Sigma $ with normal $\bar \nu$, one likes to simplify the term
\be
\mathbb{A}(V, \Sigma)  : = (\mathrm{tr}_{\bar \gamma} h)dV( \bar \nu )-V \langle \bar A,h\rangle_{\bar \gamma} .
\ee
In the asymptotically flat case, $V = 1 $ and one can take $\Sigma$ as hyperplanes which have $\bar A = 0$.
In the asymptotically hyperbolic case, it is natural to first consider those $\Sigma \subset \mathbb{H}^n$ 
which are totally umbilic, i.e.
\be \label{eq-umbilic-in-Hn}
\bar A = \frac{1}{n-1} \bar H \bar \gamma . 
\ee
Under this assumption, 
\be \label{eq-mathbbA}
\begin{split}
\mathbb{A} (V, \Sigma) 
= & \   \left( \frac{\p V}{\p \bar \nu } -\frac{\bar H}{n-1} V  \right)  \mathrm{tr}_{\bar \gamma} h \\
= & \  \left( \frac{\p }{\p \bar \nu }  \ln V -\frac{\bar H}{n-1}   \right)  V \, \mathrm{tr}_{\bar \gamma} h ,
\ \text{if} \ V \ne 0 .
\end{split}
\ee

Some basic choices of pairs $(V, \Sigma)$ making $ \mathbb{A} (V, \Sigma) = 0 $ can be described as follows. 
Consider the upper half space model of $\mathbb{H}^n$, 
$$ \mathbb{H}^n = \{ (y_1, y_2,  \ldots, y_n ) \ | \ y_1 > 0 \},   \  \hg =\frac{1}{y_1^2}(dy_1^2+\cdots +dy_n^2). $$
One knows
$$ y_1^{-1} = t - z_1  \in \mathcal{S}(\mathbb{H}^n)  . $$ 
Take $ V = y_1^{-1}$, and choose $ \Sigma$ to be 
\begin{itemize}
\item[(i)] a horosphere given by $\Sigma = \{ y_1 = c \}$ for some constant $c > 0 $; 
\item[(ii)] a totally geodesic piece that is isometric to $ \mathbb{H}^{n-1}$, say $ \Sigma = \{ y_\alpha = d \}$
for some constant $ d $, where $ \alpha \in \{ 2, \ldots, n \} $.
\end{itemize} 
It is easily seen $ \mathbb{A(}V, \Sigma) = 0 $ on such an $\Sigma$.

Making use of these facts, Jang and the author derived the following in \cite{JangMiao20}.

\begin{figure}[ht]
\centering
\def\svgwidth{3.5in}
\begingroup%
  \makeatletter%
  \providecommand\color[2][]{%
    \errmessage{(Inkscape) Color is used for the text in Inkscape, but the package 'color.sty' is not loaded}%
    \renewcommand\color[2][]{}%
  }%
  \providecommand\transparent[1]{%
    \errmessage{(Inkscape) Transparency is used (non-zero) for the text in Inkscape, but the package 'transparent.sty' is not loaded}%
    \renewcommand\transparent[1]{}%
  }%
  \providecommand\rotatebox[2]{#2}%
  \newcommand*\fsize{\dimexpr\f@size pt\relax}%
  \newcommand*\lineheight[1]{\fontsize{\fsize}{#1\fsize}\selectfont}%
  \ifx\svgwidth\undefined%
    \setlength{\unitlength}{569.88774337bp}%
    \ifx\svgscale\undefined%
      \relax%
    \else%
      \setlength{\unitlength}{\unitlength * \real{\svgscale}}%
    \fi%
  \else%
    \setlength{\unitlength}{\svgwidth}%
  \fi%
  \global\let\svgwidth\undefined%
  \global\let\svgscale\undefined%
  \makeatother%
  \begin{picture}(1,0.5622811)%
    \lineheight{1}%
    \setlength\tabcolsep{0pt}%
    \put(0,0){\includegraphics[width=\unitlength,page=1]{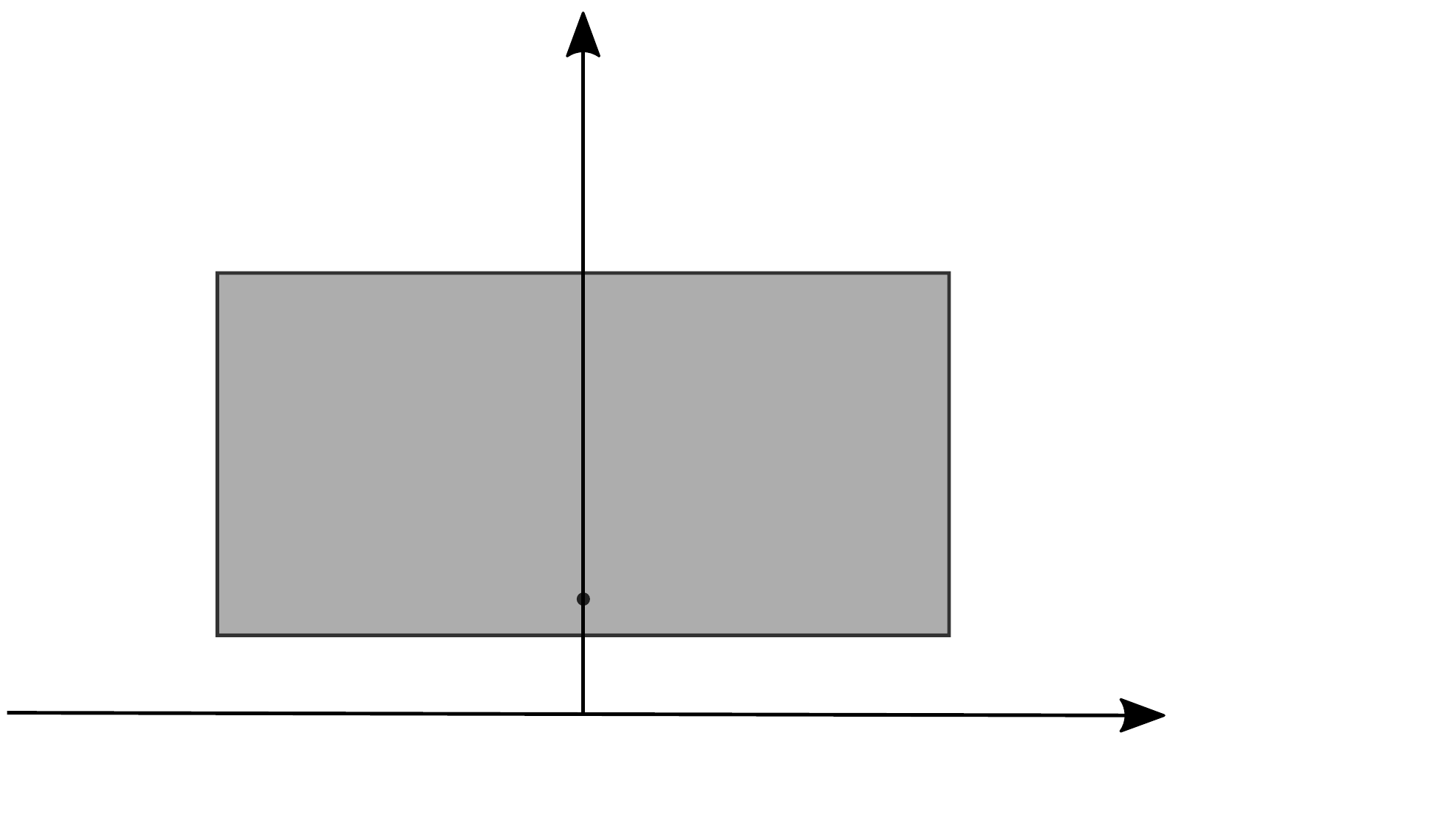}}%
    \put(0.4207261,0.53527467){\color[rgb]{0,0,0}\makebox(0,0)[lt]{\lineheight{1.25}\smash{\begin{tabular}[t]{l}$y_1$\end{tabular}}}}%
    \put(0.785564,0.02332299){\color[rgb]{0,0,0}\makebox(0,0)[lt]{\lineheight{1.25}\smash{\begin{tabular}[t]{l} 
    \end{tabular}}}}%
    \put(0.41105424,0.33015026){\color[rgb]{0,0,0}\makebox(0,0)[lt]{\lineheight{1.25}\smash{\begin{tabular}[t]{l}$e^L$\end{tabular}}}}%
    \put(0.40769948,0.0796641){\color[rgb]{0,0,0}\makebox(0,0)[lt]{\lineheight{1.25}\smash{\begin{tabular}[t]{l}$e^{-L}$\end{tabular}}}}%
    \put(0.62766768,0.00671733){\color[rgb]{0,0,0}\makebox(0,0)[lt]{\lineheight{1.25}\smash{\begin{tabular}[t]{l}$\sigma(L)$\end{tabular}}}}%
    \put(0.11929685,0.01047758){\color[rgb]{0,0,0}\makebox(0,0)[lt]{\lineheight{1.25}\smash{\begin{tabular}[t]{l}$-\sigma(L)$\end{tabular}}}}%
    \put(0,0){\includegraphics[width=\unitlength,page=2]{y-prism.pdf}}%
    \put(0.81943456,0.11726541){\color[rgb]{0,0,0}\makebox(0,0)[lt]{\lineheight{1.25}\smash{\begin{tabular}[t]{l}
    \end{tabular}}}}%
    \put(0,0){\includegraphics[width=\unitlength,page=3]{y-prism.pdf}}%
    \put(0.81112883,0.36755682){\color[rgb]{0,0,0}\makebox(0,0)[lt]{\lineheight{1.25}\smash{\begin{tabular}[t]{l}
    \end{tabular}}}}%
    \put(0,0){\includegraphics[width=\unitlength,page=4]{y-prism.pdf}}%
  \end{picture}%
\endgroup%
		\caption{A $y$-coordinate rectangular prism in $\mathbb{H}^n$}
\label{fig:y-prism}
\end{figure}

\begin{thm} [\cite{JangMiao20}] \label{thm-J-M}
Let $(M^n,g)$ be an asymptotically hyperbolic manifold.
Consider a large coordinate rectangular prism 
$$ P_L = \left\{ e^{-L} \le y_1 \le e^{L} , \ | y_\alpha | \le \sigma (L) , \  \alpha = 2, \ldots, n \right\} \subset \mathbb{H}^n . $$
Here $ \sigma (L) $ is any function on $ (L_0, \infty)$ for some $L_0 > 0 $ 
so that $ \lim_{L \to \infty}  \sigma (L) = \infty$. 
Let $ V = y_1^{-1}$.  Then
\begin{enumerate}
\item[(a)] The mass components $p_0$ and $p_1$ satisfies
\bee
p_0 - p_1 = 2 \lim_{ L \to \infty} \left[  \int_{ \mathcal{F} (\p P_L) }  V (\bar H - H) \, d  \sigma 
+  \int_{ \mathcal{E}(\p P_L  ) } V  \left( \alpha - \frac{\pi}{2} \right) \, d \mu \right]  . 
\eee
\item[(b)] If $ \sigma (L)^{n-2-2 q} = o ( e^{(q-n+1)L} ) $ as $ L \to \infty$, 
which holds automatically if $ q \ge n -1$, 
then
\bee
\lim_{ L \to \infty}  \int_{ \mathcal{E}(\p P_L  ) } V  \left( \alpha - \frac{\pi}{2} \right) \, d \mu  = 0, 
\eee
\bee
\lim_{ L \to \infty}  \int_{ \mathcal{F} (\p P_L) }  V (\bar H - H) \, d  \sigma = 
\lim_{L \to \infty} \int_{  \mathcal{F} (\p P_L) \cap  \{ y_1 = e^{-L} \}  }  V (\bar H - H) \, d  \sigma , 
\eee
and as a result, 
\bee
p_0 - p_1 = 2 \lim_{ L \to \infty}   \int_{ \ \mathcal{F} (\p P_L) \cap  \{ y_1 = e^{-L} \} }  V (\bar H - H) \, d  \sigma .
\eee
\end{enumerate}
\end{thm}

\begin{remark}
For many examples of asymptotically hyperbolic manifolds $(M^n, g)$, 
the metric fall-off rate $q$ satisfies $ q = n $. 
For instance, this is the case if $(M, g)$ is AdS-Schwarzschild, or more generally if $(M, g)$ is 
conformally compact (see \cite{XW01}). 
\end{remark}

The content of Theorem \ref{thm-J-M} was contained in \cite[Remark 3.2]{JangMiao20}.
Though its part (a) could be thought as an analogue of Theorem \ref{thm-cube-mass} 
or \ref{thm-poly-mass}, its part (b) shows that, 
in contrast to the asymptotically flat case (see Remark \ref{rem-angle-AF}),  
 the limit of the  integral of the angle deficit along the edges of $P_L$,
and the integral of the mean curvature difference on most part of the faces of $P_L$,  
do not contribute to $ p_0 - p_1$, if $ q \ge n-1$.
Exploring this fact, Jang and the author \cite{JangMiao20} 
derived an evaluation formula of $p_0 - p_1$ using horospheres $\{ y_1 = e^{-L} \}$ alone. 

It is worth of seeking other candidates $\Sigma $ on which $ \mathbb{A}(V, \Sigma) $ can be simplified
with $ V = y_1^{-1}$.
In \cite{Chai21}, Chai found that $ \mathbb{A}(V, \Sigma) = 0 $ if $\Sigma$ is any $y$-coordinate hyperplane. 
This can be seen as follows.
Given an arbitrary
$$ \Sigma^{n-1} \subset \mathbb{H}^n = \{ (y_1, y_2, \ldots, y_n) \ | \ y_1 > 0 \}, $$
if $ A_{_E} $, $ H_{_E} $ denote its second fundamental form, mean curvature
with respect to the background Euclidean metric $ d y^2  $,   
standard conformal transformation formulae give
\be \label{eq-conformal-A}
\bar A = V A_{_E} + \frac{\p V}{\p \nu_{_E} } \gamma_{_E},
\ee
where $ \nu_{_E}  $ is the Euclidean unit normal and 
$ \gamma_{_E} $ is the background metric on $ \Sigma$.
Clearly,  \eqref{eq-umbilic-in-Hn} holds if and only if 
\be \label{eq-umbilic-y}
A_{_E} = \frac{H_{_E}}{n-1} \gamma_{_E}.
\ee
Under this assumption,
\be
V^{-1} d V (\bar \nu) - \frac{\bar H}{n-1} =  -  \frac{  H_{_E}}{n-1} V^{-1} .
\ee
Therefore,  \eqref{eq-mathbbA} becomes
\be \label{eq-Axv}
\mathbb{A}(V, \Sigma) = \left\{ 
\begin{array}{ll}
0, & \text{if} \ H_{_E} = 0 \\
- \frac{H_{_E}}{n-1} \,  \tr_{\bar \gamma} h, & \text{if} \ H_{_E} \ne 0 . 
\end{array}
\right.
\ee
Here $ H_{_E}$ is necessarily a constant.
Taking $ H_{_E} = 0 $, 
one sees $ \mathbb{A}(V, \Sigma) = 0 $ along any coordinate hyperplanes $\Sigma$ in $ \{ (y_1, y_2, \ldots, y_n ) \}$.
Thus, as pointed out by Chai in \cite{Chai21}, 
an analogous statement to part (a) of Theorem \ref{thm-J-M} holds
for an exhaustion sequence consisting of  $y$-coordinate polyhedra.

\begin{remark}
One may also explore the case $ H_{_E} \neq 0 $ in \eqref{eq-Axv}, which corresponds to 
considering domains bounded by coordinate spheres in the $y$-space. 
In this case, the presence of $ \tr_{\bar \gamma} h $ in $\mathbb{A}(V, \Sigma)$ 
suggests one could proceed in a way similar to that in Proposition A.2 of \cite{JangMiao20}.
\end{remark}
 
We thus have seen that,  on asymptotically hyperbolic manifolds, there exists a connection, 
analogous to Theorem \ref{thm-cube-mass} or \ref{thm-poly-mass}, 
between the mass quantity $p_0 - p_1$ and large Riemannian polyhedra. 
However, as shown by part (b) of Theorem \ref{thm-J-M}, in many cases, the contribution of 
the angle deficit along the edges and the mean curvature difference on a portion of the faces 
is zero as the polyhedron approaches the infinity.

\begin{figure}[h]
\centering
\def\svgwidth{2.2in}
\begingroup%
  \makeatletter%
  \providecommand\color[2][]{%
    \errmessage{(Inkscape) Color is used for the text in Inkscape, but the package 'color.sty' is not loaded}%
    \renewcommand\color[2][]{}%
  }%
  \providecommand\transparent[1]{%
    \errmessage{(Inkscape) Transparency is used (non-zero) for the text in Inkscape, but the package 'transparent.sty' is not loaded}%
    \renewcommand\transparent[1]{}%
  }%
  \providecommand\rotatebox[2]{#2}%
  \newcommand*\fsize{\dimexpr\f@size pt\relax}%
  \newcommand*\lineheight[1]{\fontsize{\fsize}{#1\fsize}\selectfont}%
  \ifx\svgwidth\undefined%
    \setlength{\unitlength}{401.93589578bp}%
    \ifx\svgscale\undefined%
      \relax%
    \else%
      \setlength{\unitlength}{\unitlength * \real{\svgscale}}%
    \fi%
  \else%
    \setlength{\unitlength}{\svgwidth}%
  \fi%
  \global\let\svgwidth\undefined%
  \global\let\svgscale\undefined%
  \makeatother%
  \begin{picture}(1,0.77247734)%
    \lineheight{1}%
    \setlength\tabcolsep{0pt}%
    \put(0,0){\includegraphics[width=\unitlength,page=1]{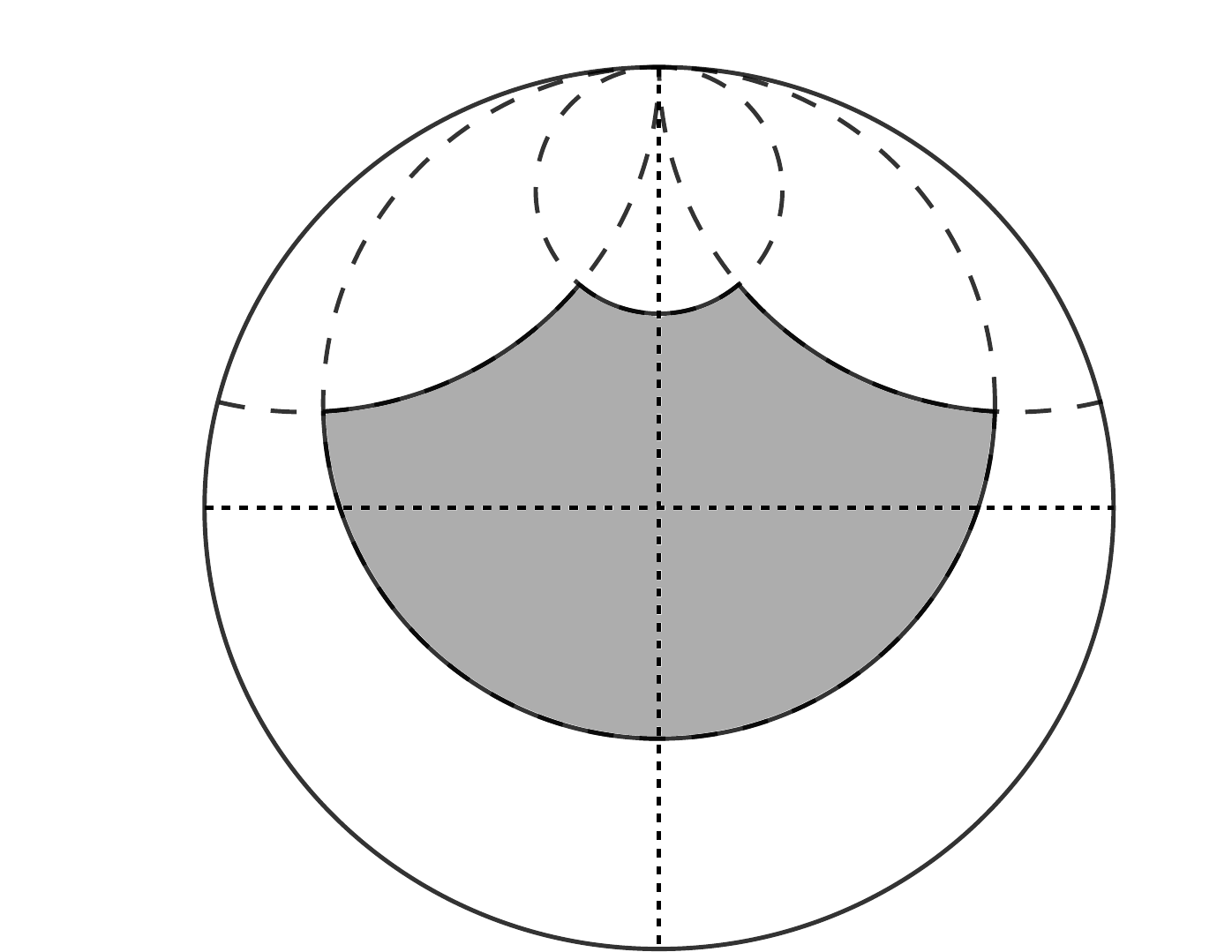}}%
    \put(0.78984457,0.71346554){\color[rgb]{0,0,0}\makebox(0,0)[lt]{\lineheight{1.25}\smash{\begin{tabular}[t]{l}
    $y_1= e^L$
    \end{tabular}}}}%
    \put(0,0){\includegraphics[width=\unitlength,page=2]{P-ball.pdf}}%
    \put(0.0147003,0.72466238){\color[rgb]{0,0,0}\makebox(0,0)[lt]{\lineheight{1.25}\smash{\begin{tabular}[t]{l}
    $y_1= e^{-L}$
    \end{tabular}}}}%
    \put(0,0){\includegraphics[width=\unitlength,page=3]{P-ball.pdf}}%
  \end{picture}%
\endgroup%
\caption{A $y$-coordinate prism in the Poincar\'{e} ball model}
\label{fig:y-prism-ball-model}
\end{figure} 

This phenomenon might be illustrated by considering the corresponding regions 
in the Poincar\'{e} ball model of ($\mathbb{H}^n, \bar g) $. 
If  $ P_L$ from Theorem \ref{thm-J-M} is illustrated in the ball model, one sees that all the edges of $P_L$, 
and all the faces of $P_L$ except that on the horosphere $\{ y_1 = e^{-L} \}$, approach a single point at the infinity. 
A quantitative relation, between the geometric convergence of the boundary of such polyhedra to the infinity
of $(M, g)$ and the analytic vanishing property of the integrals in the limit, seems worth of further exploration. 

\vspace{.3cm}

\noindent {\bf Acknowledgement}. I would like to thank Dr. Hyun Chul Jang for helpful discussions
concerning asymptotically hyperbolic manifolds.

\end{document}